# A theorem on spherically complete valued abelian groups

*Franz-Viktor Kuhlmann*

12. 6. 1997

**Abstract.** We give a criterion for a group homomorphism on a valued abelian group to be surjective and to preserve spherical completeness. We apply this to give a criterion for the existence of integration on a valued differential field. Further, we give a criterion for a sum of spherically complete subgroups of a valued abelian group to be spherically complete. This in turn can be used to determine elementary properties of power series fields in positive characteristic.

Let $(A, v)$ be a valued abelian group. That is, $v : A \ni a \mapsto va$ is a map from $A$ onto a totally ordered set with largest element $\infty$, satisfying that $va = \infty \Leftrightarrow a = 0$ and $v(a - a') \geq \min\{va, va'\}$. The latter is the **ultrametric triangle law**. We set $v(A) := \{va \mid a \in A\}$.

We recall some definitions. For $a \in A$ and $\alpha \in v(A)$, we define the **ball**

$$B_\alpha(a) := \{a' \in A \mid v(a - a') \geq \alpha\} .$$

It follows from the ultrametric triangle law that

$$B_\alpha(a) \subseteq B_{\alpha'}(a') \quad \text{if and only if} \quad a \in B_{\alpha'}(a') \text{ and } \alpha \geq \alpha' . \tag{1}$$

A set of balls in $(A, v)$ is called a **nest of balls** if it is totally ordered by inclusion. By (1), this is the case as soon as every two balls in the set have a nonempty intersection. The **intersection** of the nest is defined to be the intersection of all of its balls. Further, $(A, v)$ is called **spherically complete** if every nest of balls has a nonempty intersection.

It is well-known that the valued additive group of every power series field and, more generally, of every maximally valued field is spherically complete (cf. [P1], [P2], [PR], and [K1]). For maximally valued fields, this was essentially proved by Kaplansky in [KA], using the notion of "pseudo Cauchy sequence" instead of "nest of balls".

**Theorem:** Let $(A, v)$ and $(B, w)$ be valued abelian groups and $f : A \to B$ a group homomorphism. Assume that there is a subset $S \subseteq A$ such that:
a) the assignment

$$vs \mapsto wf(s) \qquad (s \in S)$$

defines an order preserving map $\varphi : v(S) \to w(B)$,
b) for all $a \in A$ and $s \in S$,

$$va \geq vs \Rightarrow wf(a) \geq wf(s) ,$$

c) for every $b \in B \setminus \{0\}$ there is some $s \in S$ such that

$$w(b - f(s)) > wb .$$

Suppose further that $(A, v)$ is spherically complete. Then $f$ is surjective and $(B, w)$ is spherically complete.

Proof: We will write $fa$ instead of $f(a)$. First, we observe that $w(B) \setminus \{\infty\} \subseteq \varphi(v(S))$: by assumption c), for every $b \in B$, $b \neq 0$, there is some $s \in S$ such that $w(b - fs) > wb$, which implies that $wb = wfs = \varphi vs \in \varphi(v(S))$. Since $\varphi$ is order preserving by assumption, it is also injective. Therefore, it admits an inverse $\varphi^{-1} : w(B) \setminus \{\infty\} \to v(S)$.

Take an arbitrary $a \in A$ and some $\alpha \in v(S)$; we claim that

$$f(B_\alpha(a)) \subseteq B_{\varphi\alpha}(fa) \ .$$

Indeed, if $a' \in B_\alpha(a)$ then $v(a - a') \geq \alpha$, and thus $w(fa - fa') = wf(a - a') \geq \varphi\alpha$ by assumption b), showing that $fa' \in B_{\varphi\alpha}(fa)$.

Next, using that $(A, v)$ is spherically complete we wish to show that

$$\forall a \in A \ \forall \alpha \in v(S) : \ f(B_\alpha(a)) \ = \ B_{\varphi\alpha}(fa) \tag{2}$$

where the latter is a ball in $B$. Take any $a \in A$, $\alpha \in v(S)$ and $b \in B_{\varphi\alpha}(fa)$. We have to show that $b \in f(B_\alpha(a))$. Let $B_{\alpha_i}(a_i)$, $i \in I$, be a maximal nest of balls in $B_\alpha(a)$ such that $\alpha_i \in v(S)$ and $b \in B_{\varphi\alpha_i}(fa_i)$ for all $i$. Since the nest $\{B_\alpha(a)\}$ satisfies these conditions, a maximal nest exists by Zorn's Lemma.

Since $(A, v)$ is spherically complete by assumption, there is some $\tilde{a} \in \bigcap_{i \in I} B_{\alpha_i}(a_i)$. Then $f\tilde{a} \in f(B_{\alpha_i}(a_i)) \subseteq B_{\varphi\alpha_i}(fa_i)$ and thus $w(b - f\tilde{a}) \geq \varphi\alpha_i$ for all $i$. Suppose that $b \neq f\tilde{a}$. Then by assumption c), there is some $s \in S$ such that $w(b - f\tilde{a}) < w((b - f\tilde{a}) - fs) = w(b - f(\tilde{a} + s)) =: \beta'$. We set $a' := \tilde{a} + s$ and observe that $vs = \varphi^{-1} wfs = \varphi^{-1} w(b - f\tilde{a}) \geq \alpha_i$ for all $i$. Hence, $a' = \tilde{a} + s \in B_{\alpha_i}(\tilde{a}) = B_{\alpha_i}(a_i)$ for all $i$.

If $b = fa'$ then $b \in f(B_\alpha(a))$ and we are done. Otherwise $\beta' \neq \infty$, and we can set $\alpha' := \varphi^{-1} \beta'$. As $\beta' > w(b - f\tilde{a}) \geq \varphi\alpha_i$, we have that $\alpha' > \alpha_i$ for all $i$. This shows that $B_{\alpha'}(a')$ is properly contained in all $B_{\alpha_i}(a_i)$. On the other hand, $b \in B_{\beta'}(fa') = B_{\varphi\alpha'}(fa')$. Therefore, we can enlarge our nest of balls by adding $B_{\alpha'}(a')$. As this contradicts the maximality, we conclude that $b = f\tilde{a} \in f(B_\alpha(a))$, as desired.

To show the surjectivity of $f$, take any $b \in B$. If $b = 0$ then $b = f0 \in f(A)$. So assume that $b \neq 0$. Then by assumption c), we can choose some $a \in S$ such that $\beta := w(b - fa) > wb$. If $fa = b$ then we are done. Otherwise $\beta \neq \infty$ and we can set $\alpha := \varphi^{-1}\beta \in v(S)$, obtaining that $b \in B_\beta(fa) = f(B_\alpha(a))$. This shows the surjectivity. Let us note that $\beta > wb = wfa$ and hence, $\alpha = \varphi^{-1}\beta > \varphi^{-1} wfa = va$. So we have proved, for all $b \in B$:

$$b \in f(S \cup \{0\}) \ \lor \ \exists a \in S \ \exists \alpha \in v(S) : \ \alpha > va \ \land \ b \in f(B_\alpha(a)) \ . \tag{3}$$

Finally, let us prove that $(B, w)$ is spherically complete. Take any nest of balls $\{B_{\beta_j}(b_j) \mid j \in J\}$ in $(B, w)$. We have to show that this nest has a nonempty intersection. We may assume that $\beta_j \neq \infty$ for all $j \in J$; indeed, if $\beta_{j_0} = \infty$ for some $j_0 \in J$, then the intersection is $\{b_{j_0}\} \neq \emptyset$ and there is nothing more to show.

We claim that there exists a nest of balls $\{B_{\alpha_i}(a_i) \mid i \in I\}$ in $(A, v)$ maximal with the property that for every $i \in I$,

$$\alpha_i \in v(S) \ \text{ and there is some } j \in J \text{ such that } B_{\beta_j}(b_j) \subseteq f(B_{\alpha_i}(a_i)) \ .$$

Indeed, we can choose an arbitrary $j \in J$, set $\alpha := \varphi^{-1}\beta_j$ and then, by the surjectivity of $f$, choose some $a \in A$ such that $fa = b_j$. Then by (2), $B_{\beta_j}(b_j) = f(B_\alpha(a))$. So $\{B_\alpha(a)\}$ is a nest with the required property. Hence, a maximal nest $\{B_{\alpha_i}(a_i) \mid i \in I\}$ with this property exists by Zorn's Lemma. It follows that

$$\bigcap_{j \in J} B_{\beta_j}(b_j) \subseteq \bigcap_{i \in I} f(B_{\alpha_i}(a_i)) = \bigcap_{i \in I} B_{\varphi\alpha_i}(fa_i) \,. \tag{4}$$

We wish to show that

$$f\left(\bigcap_{i \in I} B_{\alpha_i}(a_i)\right) \subseteq \bigcap_{j \in J} B_{\beta_j}(b_j) \,. \tag{5}$$

Suppose the contrary. Since $f\left(\bigcap_{i \in I} B_{\alpha_i}(a_i)\right) \subseteq \bigcap_{i \in I} f(B_{\alpha_i}(a_i))$, we must then have that inequality holds in (4) and that there is an element $\tilde{a} \in \bigcap_{i \in I} B_{\alpha_i}(a_i)$ and some $j_0 \in J$ such that $f\tilde{a} \notin B_{\beta_{j_0}}(b_{j_0}) \subsetneq \bigcap_{i \in I} B_{\varphi\alpha_i}(fa_i)$. The inequality implies that $\beta_{j_0} > \varphi\alpha_i$ for all $i$. As $b_{j_0}, f\tilde{a} \in B_{\varphi\alpha_i}(fa_i)$, we have that $w(b_{j_0} - f\tilde{a}) \geq \varphi\alpha_i$ for all $i$. By assumption c), there is some $s \in S$ such that $w(b_{j_0} - f\tilde{a}) < w((b_{j_0} - f\tilde{a}) - fs) = w(b_{j_0} - f(\tilde{a}+s)) =: \tilde{\beta}$. As before, we set $a' := \tilde{a} + s$ and observe that $a' \in B_{\alpha_i}(a_i)$ for all $i$. We set $\beta' := \min\{\tilde{\beta}, \beta_{j_0}\} < \infty$ and $\alpha' := \varphi^{-1}\beta'$. Then $\beta' > \varphi\alpha_i$ and $\alpha' > \alpha_i$ for all $i$. It follows that $B_{\alpha'}(a')$ is properly contained in all $B_{\alpha_i}(a_i)$ for all $i$. On the other hand, $\beta_{j_0} \geq \beta'$ and $w(b_{j_0} - fa') = \tilde{\beta} \geq \beta'$ imply that $B_{\beta_{j_0}}(b_{j_0}) \subseteq B_{\beta'}(fa') = f(B_{\alpha'}(a'))$. Therefore, we can enlarge our nest of balls by adding $B_{\alpha'}(a')$. As this contradicts the maximality, we conclude that (5) must hold.

Since $(A,v)$ is spherically complete by assumption, $\bigcap_{i \in I} B_{\alpha_i}(a_i)$ is nonempty. Hence by (5), $\bigcap_{j \in J} B_{\beta_j}(b_j)$ is nonempty, as desired. □

There is also a version of this theorem for ultrametric spaces, where addition is not available ([K3]).

**Remarks:** Assume that the conditions of the theorem hold. Then in particular, (3) holds for every $b \in B$.

Take any $a \in S$ and $\alpha \in v(S)$ such that $va < \alpha$. It follows that every $a' \in B_\alpha(a)$ has value $va' = \min\{va, \alpha\} = va$. By (2), $fa' \in f(B_\alpha(a)) = B_{\varphi\alpha}(fa)$. Since $\varphi\alpha > wfa$ and $a \in S$, we have that $wfa' = \min\{wfa, w(fa' - fa)\} = wfa = \varphi va = \varphi va'$. Hence, the subset $S \cup B_\alpha(a) \subseteq A$ still has property a). Trivially, it also has properties b) and c). Therefore, every maximal subset $S \subseteq A$ having properties a), b) and c) will contain every ball $B_\alpha(a)$ of the above form. By (3), it will therefore satisfy that

$$f(S \cup \{0\}) \;=\; B \,.$$

On the other hand, if $S$ can be chosen such that $S - S \subseteq S$ (which in fact means that $S$ is a subgroup), then $f : S \to B$ will be injective. Indeed, if $s, s' \in S$ such that $fs = fs'$, then $f(s - s') = 0$ and $wf(s - s') = \infty$. With $s - s' \in S$ and $0 \in S$, this implies that $v(s - s') = \infty$, i.e., $s = s'$. If in addition for all $a \in A$ there is some $s \in S$ such that $fs = fa$, then $f : S \to B$ will be an isomorphism of valued groups. This is the situation of application (A) below.

Note that if $S$ contains an element $s_0$ such that $fs_0 = 0$ (e.g., if $0 \in S$), then by the first part of our proof, $\varphi : v(S) \to w(B)$ will be an order preserving bijection. For this, it is not necessary that $(A, v)$ be spherically complete.

**Applications:**
(A) Let $(K, D)$ be a differential field with field of constants $C = \{a \in K \mid Da = 0\}$. Following M. Rosenlicht [R1], a valuation $v$ of $K$ is called a **differential valuation** if it is trivial on $C$ and satisfies

$$\forall a, b \in K : va \geq 0 \land vb > 0 \land b \neq 0 \Rightarrow v\left(\frac{bDa}{Db}\right) > 0 . \tag{6}$$

These conditions imply that

$$\forall a, b \in K \setminus \{0\}, va \neq 0, vb \neq 0 : va \leq vb \Leftrightarrow vDa \leq vDb . \tag{7}$$

If $C$ is a set of representatives for the $v$-residues of $K$, then also the converse holds, i.e., $v$ is differential if and only if it satisfies (7), cf. [R1].

We say that $(K, D)$ **admits integration** if $D$ is surjective. If $v$ is a differential valuation on $K$, then following Rosenlicht [R2], we say that $(K, D)$ **admits asymptotic integration** if for every $b \in K$, $b \neq 0$, there is some $a \in K$ such that $v(b - Da) > vb$.

Let $k$ be any field and $G$ any ordered abelian group. The **(generalized) power series field** $k((G))$ consists of all formal sums $a = \sum_{g \in G} c_g t^g$ with $c_g \in k$ and well-ordered **support** $\text{supp}(a) = \{g \in G \mid c_g \neq 0\}$. The **canonical valuation** on $k((G))$ is given by $va := \min \text{supp}(a) \in G$ and $v0 := \infty$.

The following is applied in [KK] to prove integration for exponential-logarithmic power series:

*Let $k((G))$ be a power series field, endowed with a derivation $D$ for which $k$ is the field of constants. Suppose that the canonical valuation $v$ of $K$ is a differential valuation (with respect to $D$). Then $(k((G)), D)$ admits integration if and only if it admits asymptotic integration.*

Proof: As we have already mentioned, the additive group of the power series field $k((G))$, endowed with the canonical valuation $v$, is spherically complete. We let $S$ be the subgroup of all power series $a$ in $k((G))$ without constant term (that is, if $a = \sum_{g \in G} c_g t^g$ then $c_0 = 0$). Then $k((G)) = S \oplus k$ and $v(S) = G \setminus \{0\}$. Therefore, (7) implies assumption a) of the theorem.

If $a = s' + c \in k((G))$ with $s' \in S$ and $c \in k$, then $va = vs' < 0$ or $vs' \geq va \geq 0$ (with $vs' = va > 0$ if $c = 0$). In all cases, $vs' \geq va$. Hence if $s \in S$ such that $va \geq vs$, then $vs' \geq vs$ and by (7), $vDa = vDs' \geq vDs$. This shows that also assumption b) holds.

Since for every $a \in k((G))$ there is some $s \in S$ such that $Da = Ds$, assumption c) holds if and only if $(k((G)), D)$ admits asymptotic integration. Now the theorem proves one direction, and the other one is clear since integration implies asymptotic integration. □

In this application, $S$ is a subgroup of the additive group of $k((G))$, having the property that for every $a \in k((G))$ there is some $s \in S$ such that $Ds = Da$. Hence by the remark above, $D : S \to k((G))$ is an isomorphism of valued abelian groups. In particular, this shows that for every $b \in k((G))$ there is a *unique* $s \in S$ such that $Ds = b$.

On the other hand, $S$ is not a maximal subset of $k((G))$ having properties a), b) and c). Indeed, as $v(S) = G \setminus \{0\}$, we can pick $a \in S$ and $\alpha \in v(S)$ such that $va < \alpha < 0$. Then $a + c \in B_\alpha(a)$ for every $c \in k$ and thus, $B_\alpha(a) \not\subseteq S$. However, $S$ is a maximal *subgroup* having properties a), b) and c).

(B) Let $(\mathcal{A}, v)$ be a valued abelian group and $A_1, \ldots, A_n$ be subgroups of $\mathcal{A}$. We call the sum $A_1 + \ldots + A_n \subseteq \mathcal{A}$ **pseudo-direct** if for every $a \in A_1 + \ldots + A_n$, $a \neq 0$, there are $a_i \in A_i$ such that

$$v \sum_{i=1}^n a_i = \min_{1 \leq i \leq n} va_i \quad \text{and} \quad v\left(a - \sum_{i=1}^n a_i\right) > va \;.$$

The following is applied in [K1], [K2] to determine elementary properties of the power series field $\mathbb{F}_p((t))$:

*Assume that the subgroups $(A_i, v)$ of $(\mathcal{A}, v)$, $1 \leq i \leq n$, are spherically complete. If the sum $A_1 + \ldots + A_n \subseteq \mathcal{A}$ is pseudo-direct, then it is also spherically complete.*

Proof: We endow the direct product $A := A_1 \times \ldots \times A_n$ with the **minimum valuation** by setting

$$v(a_1, \ldots, a_n) := \min_{1 \leq i \leq n} va_i \;.$$

Further, we set $B := A_1 + \ldots + A_n \subseteq \mathcal{A}$ and take $f : A \to B$ to be the group homomorphism defined by $f(a_1, \ldots, a_n) := \sum_{i=1}^n a_i$. We let $S \subseteq A$ be the set of all $n$-tuples $(a_1, \ldots, a_n)$ which satisfy that $v \sum_{i=1}^n a_i = \min_{1 \leq i \leq n} va_i$. We note that $v(S) = v(A) \subseteq v(\mathcal{A})$ and that the map $\varphi$ defined as in assumption a) of the theorem is simply the identity on $v(S)$. Assumption b) of the theorem holds since by the ultrametric triangle law, we always have that $v \sum_{i=1}^n a_i \geq \min_{1 \leq i \leq n} va_i$. Finally, the sum $A_1 + \ldots + A_n \subseteq \mathcal{A}$ is pseudo-direct if and only if assumption c) of the theorem holds. As it is well-known that the direct product of finitely many spherically complete valued abelian groups, endowed with the minimum valuation, is again spherically complete (cf. [K1]), our assertion now follows from the theorem. $\square$

In this application, $S$ will in general not be a subgroup of $A$, and $f$ will in general not be injective. Indeed, the groups $A_i$ may have non-trivial intersection, so two elements of $S$ may be equal up to a permutation of the coördinates $a_i$.

(C) Let $(A, v)$ be a spherically complete valued abelian group, and $\alpha \in v(A)$. We leave it to the reader to define a canonical valuation on the quotient group $A/B_\alpha(0)$ and to apply the theorem in order to show that with this valuation, $A/B_\alpha(0)$ is spherically complete.

# References

[KA]  Kaplansky, I.: *Maximal fields with valuations I*, Duke Math. J. **9** (1942), 303–321

[K1]  Kuhlmann, F.-V.: *Valuation theory of fields, abelian groups and modules*, preprint, Heidelberg (1996), to appear in the "Algebra, Logic and Applications" series (Gordon and Breach), eds. A. Macintyre and R. Göbel

[K2]  Kuhlmann, F.-V.: *Elementary properties of power series fields over finite fields*, conference hand-out, Toronto (January 1997); prepublication in: Structures Algébriques Ordonnées, Séminaire Paris VII (1997), and in: The Fields Institute Preprint Series (June 1997)

[K3]  Kuhlmann, F.-V.: *Ultrametric spaces*, in preparation

[KK]  Kuhlmann, F.-V. – Kuhlmann, S.: *Exponential-logarithmic power series*, in preparation

[P1]  Prieß-Crampe, S.: *Angeordnete Strukturen. Gruppen, Körper, projektive Ebenen*, Ergebnisse der Mathematik und ihrer Grenzgebiete **98**, Springer (1983)

[P2]  Prieß-Crampe, S.: *Der Banachsche Fixpunktsatz für ultrametrische Räume*, Results in Mathematics **18** (1990), 178–186

[PR]  Prieß-Crampe, S. – Ribenboim, P.: *Fixed Points, Combs and Generalized Power Series*, Abh. Math. Sem. Hamburg **63** (1993), 227-244

[R1]  Rosenlicht, M.: *Differential valuations*, Pacific J. Math. **86** (1980), 301–319

[R2]  Rosenlicht, M.: *On the value group of a differential valuation*, Amer. J. Math. **191** (1979), 258–266

The Fields Institute, 222 College Street, Toronto, Ontario M5T 3J1, Canada
email: fkuhlman@fields.utoronto.ca